\documentclass[reqno,11pt]{amsart}
\usepackage{citesort}

\newtheorem{theorem}{Theorem}[section]
\newtheorem{lemma}[theorem]{Lemma}

\theoremstyle{definition}
\newtheorem{definition}[theorem]{Definition}
\newtheorem{notation}[theorem]{Notation}

\theoremstyle{remark}

\numberwithin{equation}{section}

\newcommand{\abs}[1]{\lvert#1\rvert}
\def\R{{\mathbb{R}}} 
\def\N{{\mathbb{N}}}
\def\epsilon{{\varepsilon}}
\def\phi{{\varphi}}
\def\theta{{\vartheta}}
\DeclareMathOperator{\tr}{tr}
\DeclareMathOperator{\dist}{dist}
\def\trw{{\tr w^{ij}}} 
\def\fracp#1#2{\frac{\partial #1}{\partial #2}}
\def\ol#1{{\overline{#1}}}
\def\ul#1{{\underline{#1}}}

\begin{document}

\title{Schouten Tensor Equations in Conformal Geometry 
  with Prescribed Boundary Metric}

\author{Oliver C. Schn\"urer}
\address{Max Planck Institute for Mathematics in the Sciences, 
  Inselstr.\ 22-26, 04103 Leipzig, Germany}
\curraddr{}
\email{Oliver.Schnuerer@mis.mpg.de}


\date{July 2003.}


\begin{abstract}
On a manifold with boundary, we deform the metric 
conformally. 
This induces a deformation of the Schouten tensor. 
We fix the metric at the boundary and realize
a prescribed value for the product of the eigenvalues
of the Schouten tensor in the interior, 
provided that there exists a subsolution.
\end{abstract}

\maketitle

\markboth{OLIVER C. SCHN\"URER}{SCHOUTEN EQUATIONS WITH BOUNDARY}

\section{Introduction}

Let $\left(M^n,g_{ij}\right)$ be an $n$-dimensional Riemannian 
manifold, $n\ge3$.
The Schouten tensor $(S_{ij})$ of $(M^n,g_{ij})$ is defined as
$$S_{ij}=\tfrac1{n-2}\left(R_{ij}-\tfrac1{2(n-1)}Rg_{ij}\right),$$
where $(R_{ij})$ and $R$ denote the Ricci and scalar curvature
of $(M^n,g_{ij})$, respectively. 
Consider the manifold $\left(\tilde M^n,\,\tilde g_{ij}\right)
=\left(M^n,\,e^{-2u}g_{ij}\right)$, where we have used 
$u\in C^2(M^n)$ to deform the metric conformally.
The Schouten tensors $S_{ij}$ of $g_{ij}$ and 
$\tilde S_{ij}$ of $\tilde g_{ij}$ are related by
$$\tilde S_{ij}=u_{ij}+u_iu_j-\tfrac12\vert\nabla u\vert^2
g_{ij}+S_{ij},$$
where indices of $u$ denote covariant derivatives
with respect to the background metric $g_{ij}$, moreover
$\vert\nabla u\vert^2=g^{ij}u_iu_j$ and $\left(g^{ij}\right)
=(g_{ij})^{-1}$. Eigenvalues of the Schouten tensor are 
computed with respect to the background metric $g_{ij}$, so
the product of the eigenvalues of the Schouten tensor 
$(\tilde S_{ij})$ equals a given function $s:M^n\to\R$, if
\begin{equation}\label{s eqn}
\frac{\displaystyle\det\left(u_{ij}+u_iu_j-\tfrac12\vert\nabla u\vert^2
g_{ij}+S_{ij}\right)}{\displaystyle e^{-2nu}\det\left(g_{ij}\right)}=s(x).
\end{equation}
We say that $u$ is an admissible solution for \eqref{s eqn}, if
the tensor in the determinant in the numerator is positive 
definite. At admissible solutions, \eqref{s eqn} becomes an
elliptic equation.
As we are only interested in admissible solutions,
we will always assume that $s$ is positive.

Let now $M^n$ be compact with boundary and 
$\ul u:M^n\to\R$ be a smooth (up to the boundary) admissible 
subsolution to \eqref{s eqn}
\begin{equation}\label{sub sol}
\frac{\displaystyle\det\left(\ul u_{ij}+\ul u_i\ul u_j
-\tfrac12\vert\nabla \ul u\vert^2 g_{ij}+S_{ij}\right)}
{\displaystyle e^{-2n\ul u}\det\left(g_{ij}\right)}\ge s(x).
\end{equation}

Assume that there exists a
supersolution $\ol u$ to \eqref{s eqn} 
fulfilling some technical conditions specified in
Definition \ref{sup def}.
Assume furthermore that $M^n$ admits a strictly 
convex function $\chi$. Without loss of generality,
we have $\chi_{ij}\ge g_{ij}$ for the second covariant
derivatives of $\chi$ in the matrix sense. 
 
The conditions of the preceding paragraph are
automatically fulfilled if $M^n$ is a compact subset 
of flat $\R^n$ and $\ul u$ fulfills
\eqref{sub sol} and in addition
$\det(\ul u_{ij})\ge s(x)e^{-2n\ul u}\det(g_{ij})$
with $\ul u_{ij}>0$ in the matrix sense.
Then Lemma \ref{sup sol exi lem}
implies the existence of a supersolution
and we may take $\chi=\vert x\vert^2$.

We impose the boundary condition that the metric
$\tilde g_{ij}$ at the boundary is prescribed,
$$\tilde g_{ij}=e^{-2\ul u}g_{ij}\quad\text{on~}
\partial M^n.$$ 

Assume that all data are smooth up to the boundary.
We prove the following
\begin{theorem}
Let $M^n$, $g_{ij}$, $\ul u$, $\ol u$, $\chi$, and $s$ be as above. 
Then there exists a metric $\tilde g_{ij}$, 
conformally equivalent to $g_{ij}$ with $\tilde g_{ij}
=e^{-2\ul u}g_{ij}$ on $\partial M^n$ such that 
the product of the eigenvalues of the Schouten
tensor induced by $\tilde g_{ij}$ equals $s$.
\end{theorem}
This follows readily from
\begin{theorem}\label{main thm}
Under the assumptions stated above, there exists 
an admissible function $u\in C^0(M^n)\cap 
C^\infty(M^n\setminus\partial M^n)$ solving 
\eqref{s eqn} such that $u=\ul u$ on $\partial M^n$.
\end{theorem}

Recently, in a series of papers, Jeff Viaclovsky studied 
conformal deformations of metrics on closed manifolds 
and elementary symmetric functions $S_k$, $1\le k\le n$,
 of the eigenvalues
of the associated Schouten tensor, see e.\ g.\ 
\cite{ViaclovskyCAG} for existence results. Pengfei Guan,
Jeff Viaclovsky, and Guofang Wang provide an estimate
that can be used to show compactness of manifolds 
with lower bounds on elementary symmetric functions of
the eigenvalues of the Schouten tensor 
\cite{GVWSchoutenCompact}. 
A similar equation arises in
geometric optics \cite{WangInverse,GuanWang}. 
Xu-Jia Wang proved the existence of solutions to 
Dirichlet boundary value problems for an equation 
similar to \eqref{s eqn}, 
provided that the domains are small.
In \cite{OSRefl} we provide a transformation that shows
the similarity between reflector and Schouten tensor
equations. Pengfei Guan
and Xu-Jia Wang obtained local $C^2$-estimates 
\cite{GuanWang}.
This was extended by Pengfei Guan and Guofang Wang
to local $C^1$- and $C^2$-estimates in the case
of elementary
symmetric functions $S_k$ of the Schouten tensor
of a conformally deformed metric \cite{GuanLocalSchouten}. 
Boundary value problems for Monge-Amp\`ere equations
have been studied by Luis Caffarelli, Louis Nirenberg,
and Joel Spruck in \cite{CNS1} an many other people later 
on. For us, those articles using subsolutions as used
by Bo Guan and Joel Spruck will be especially useful
\cite{GuanSpruck,BGuanTrans,OSMathZ,NehringCrelle}.

It follows directly from the proof that we can also 
solve Equation \eqref{s eqn} on non-compact complete 
manifolds provided that there exist appropriate sub- and
supersolutions with locally bounded difference in $C^0$.
Then we can solve \eqref{s eqn} with an artificially
introduced Dirichlet boundary condition on a sequence
of growing domains exhausting the non-compact manifold.
A subsequence of these solutions converges then to a 
solution of \eqref{s eqn} on the manifold.
This works as the local $C^2$-estimates
in \cite{GuanLocalSchouten} depend only on 
a local bound for $\abs u$.
Note that either $s(x)$ has to decay at
infinity or the manifold with metric $e^{-2u}g_{ij}$
is non-complete. Otherwise, \cite{GVWSchoutenCompact}
implies a positive lower bound on the Ricci tensor,
i.\ e.\ $\tilde R_{ij}\ge\frac1c\tilde g_{ij}$ 
for some positive
constant $c$. This yields compactness of the manifold
\cite{GroKliMey}.
 
It is a further issue to solve similar problems
for other elementary symmetric functions of the
Schouten tensor. As the induced mean curvature
of $\partial M^n$ is related to the Neumann 
boundary condition, this is another natural 
boundary condition.

To show existence for a boundary value problem for
fully nonlinear equations like Equation 
\eqref{s eqn}, one usually proves $C^2$-estimates 
up to the boundary. Then standard results imply 
$C^k$-bounds for $k\in\N$ and existence results.
In our situation, however, we don't expect that
$C^2$-estimates up to the boundary can be proved.
This is due to the gradient terms appearing in
the determinant in \eqref{s eqn}. It is possible
to overcome these difficulties by considering only 
small domains \cite{WangInverse}. Our method is
different. We regularize the equation and prove
full regularity up to the boundary
for the regularized equation. Then we use
the fact, that interior $C^k$-estimates \cite{GuanLocalSchouten}
can be obtained independent of the regularization.
Moreover, we can prove uniform $C^1$-estimates.
Thus we can pass to a limit and get a solution
in $C^0(M^n)\cap C^\infty(M^n\setminus\partial M^n)$.

To be more precise, we rewrite \eqref{s eqn} in the
form
\begin{equation}\label{f eqn}
\log\det\left(u_{ij}+u_iu_j-\tfrac12\vert\nabla u
\vert^2g_{ij}+S_{ij}\right)=f(x,u),
\end{equation}
where $f\in C^\infty(M^n\times\R)$. Our method
can actually be applied to any equation of that form
provided that we have sub- and supersolutions. Thus
we consider in the following equations of the 
form \eqref{f eqn}.
Equation \eqref{f eqn} makes sense in any dimension
provided that we replace $S_{ij}$ by a smooth
tensor. In this case Theorem \ref{main thm} is valid
in any dimension. Note that even without the factor 
$\frac1{n-2}$ in the definition of the Schouten tensor,
our equation is not elliptic for $n=2$ for any function
$u$ as the trace $g^{ij}(R_{ij}-\frac12Rg_{ij})$ equals
zero, so there has to be a non-positive eigenvalue 
of that tensor.
Let $\psi:M^n\to[0,1]$ be smooth, $\psi=0$ in a
neighborhood of the boundary. Then our strategy 
is as follows. We 
consider a sequence $\psi_k$ of those functions 
that fulfill $\psi_k(x)=1$ for 
$\dist(x,\partial M^n)>\tfrac2k$, $k\in\N$, and
boundary value problems
\begin{equation}\label{psi eqn}
\left\{\begin{array}{rcll}
\log\det\left(u_{ij}+\psi u_iu_j-\tfrac12\psi
\vert\nabla u\vert^2g_{ij}+T_{ij}\right)&=&
f(x,u)&\mbox{in~}M^n,\\
u&=&\ul u&\mbox{on~}\partial M^n.
\end{array}\right.
\end{equation}
We dropped the index $k$ to keep the notation simple.
The tensor $T_{ij}$ coincides with $S_{ij}$ on
$\left\{x\in M^n:\dist(x,\partial M^n)>\tfrac2k\right\}$
and interpolates smoothly to $S_{ij}$ plus a sufficiently large 
constant multiple of the background metric $g_{ij}$
near the boundary.

Our sub- and supersolutions act as barriers and  
imply uniform $C^0$-estimates.
We prove uniform $C^1$-estimates based on the 
admissibility of solutions. Admissibility means here
that $u_{ij}+\psi u_iu_j-\tfrac12\psi
\vert\nabla u\vert^2+S_{ij}$ is positive definite
for those solutions. As mentioned above, we can't
prove uniform $C^2$-estimates for $u$, but we get
$C^2$-estimates that depend on $\psi$. These estimates
guarantee, that we can apply standard methods
(Evans-Krylov-Safanov theory, Schauder estimates
for higher derivatives, and mapping degree theory
for existence, see e.\ g.\ \cite{GT,Taylor3,BGuanTrans})
to prove existence of a smooth admissible
solution to \eqref{psi eqn}.
Then we use \cite{GuanLocalSchouten} to get uniform
interior a priori estimates on compact subdomains
of $M^n$ as $\psi=1$ in a neighborhood of these 
subdomains for all but a finite number of regularizations.
These a priori estimates suffice to pass to a subsequence
and to obtain an admissible solution to \eqref{f eqn}
in $M^n\setminus\partial M^n$. As $u^k=u=\ul u$
for all solutions $u^k$ of the regularized equation and 
those solutions have uniformly bounded gradients, the
boundary condition is preserved when we pass to the
limit and we obtain Theorem \ref{main thm} provided that
we can prove $\left\Vert u^k\right\Vert_{C^1(M^n)}\le c$ 
uniformly and $\left\Vert u^k\right\Vert_{C^2(M^n)}
\le c(\psi)$. These estimates are proved in Lemmata
\ref{unif C1} and \ref{int C2}, the crux of this paper.

The rest of the article is organized as follows.
We introduce supersolutions and some notation in
Section \ref{upp barr}. We mention $C^0$-estimates
in Section \ref{C0 sec}.
In Section \ref{C1 sec}, we prove uniform 
$C^1$-estimates. 
Then the $C^2$-estimates proved in
Section \ref{C2 sec} complete the a priori estimates
and the proof of Theorem \ref{main thm}.

The author wants to thank J\"urgen Jost and the 
Max Planck Institute for Mathematics in the 
Sciences for support and Guofang Wang for interesting
discussions about the Schouten tensor.


\section{Supersolutions and Notation}
\label{upp barr}

Before we define a supersolution, we explain
more explicitly, how we regularize the equation. For 
fixed $k\in\N$ we take $\psi_k$ such that
$$\psi_k(x)=\begin{cases}
0 & \dist(x,\partial M^n)<\tfrac 1k,\\
1 & \dist(x,\partial M^n)>\tfrac 2k
\end{cases}$$
and $\psi_k$ is smooth with values in $[0,1]$.
Again, we drop the index $k$ to keep the notation simple.
We fix $\lambda\ge0$ sufficiently large so that
\begin{equation}\label{mod sub sol}
\log\det\left(\ul u_{ij}+\psi\ul u_i\ul u_j-\tfrac12
\psi\vert\nabla\ul u\vert^2g_{ij}+S_{ij}+\lambda(1-\psi)
g_{ij}\right)\ge f(x,\ul u)
\end{equation}
for any $\psi=\psi_k$, independent
of $k$. 
As $\log\det(\cdot)$ is a concave function on positive
definite matrices, \eqref{mod sub sol} follows for
$k$ sufficiently large, if
$$\log\det\left(\ul u_{ij}+\ul u_i\ul u_j-\tfrac12
\vert\nabla u\vert^2g_{ij}+S_{ij}\right)\ge f(x,\ul u)
\quad\text{on~}M^n$$
and 
$$\log\det\left(\ul u_{ij}+S_{ij}+\lambda g_{ij}\right)
\ge f(x,\ul u)\quad\text{near~}\partial M^n,$$
provided that the arguments of the determinants are
positive definite.

We define
\begin{definition}[supersolution]\label{sup def}
A smooth function $\ol u:M^n\to\R$ is called a 
supersolution, if $\ol u\ge\ul u$ and for any
$\psi$ as considered above, 
$$\log\det\left(\ol u_{ij}+\psi\ol u_i\ol u_j-\tfrac12
\psi\vert\nabla\ol u\vert^2g_{ij}+S_{ij}+\lambda(1-\psi)
g_{ij}\right)\le f(x,\ul u)$$
holds for those points in $M^n$ for which the tensor in 
the determinant is positive definite. 

\end{definition}
\begin{lemma}\label{sup sol exi lem}
If $M^n$ is a compact subdomain of flat $\R^n$,
the subsolution $\ul u$ fulfills \eqref{sub sol}
and in addition
$$\det(\ul u_{ij})\ge s(x)e^{-2n\ul u}
\det(g_{ij})$$
holds, where $\ul u_{ij}>0$ in the matrix sense, 
then there exists a supersolution.
\end{lemma}
\begin{proof}
In flat $\R^n$, we have $S_{ij}=0$. The inequality
\begin{equation}\label{middle}
\frac{\displaystyle\det\left(\ul u_{ij}+\psi\ul u_i\ul u_j
-\tfrac12\psi\vert\nabla \ul u\vert^2 g_{ij}\right)}
{\displaystyle e^{-2n\ul u}\det\left(g_{ij}\right)}\ge s(x)
\end{equation}
is fulfilled if $\psi$ equals $0$ or $1$ by assumption.
As above, \eqref{middle} follows
for any $\psi\in[0,1]$. Thus \eqref{mod sub sol} is
fulfilled for $\lambda=0$.

Let $\ol u=\sup\limits_{M^n}\ul u+1+\epsilon\vert x\vert^2$
for $\epsilon>0$. It can be verified directly that $\ol u$
is a supersolution for $\epsilon>0$ fixed sufficiently
small.
\end{proof}

\begin{notation}
We set 
\begin{align*}
w_{ij}=&u_{ij}+\psi u_iu_j-\tfrac12\psi\vert\nabla u\vert^2
g_{ij}+S_{ij}+\lambda(1-\psi) g_{ij}\\
=&u_{ij}+\psi u_iu_j-\tfrac12\psi\vert\nabla u\vert^2
g_{ij}+T_{ij}
\end{align*}
and use $\left(w^{ij}\right)$ to denote the inverse of $(w_{ij})$.
The Einstein summation convention is used. We lift and lower
indices using the background metric. Vectors of length one
are called directions. Indices, sometimes 
preceded by a semi-colon, denote covariant derivatives.
We use indices preceded by a colon for partial derivatives.
Christoffel symbols of the background metric are denoted
by $\Gamma^k_{ij}$, so $u_{ij}=u_{;ij}=u_{,ij}-\Gamma^k_{ij}u_k$.
Using the Riemannian curvature tensor $(R_{ijkl})$, 
we can interchange covariant differentiation
\begin{equation}\label{interchange}
\begin{split}
u_{ijk}=&u_{kij}+u_a g^{ab}R_{bijk},\\
u_{iklj}=&u_{ikjl}+u_{ka}g^{ab}R_{bilj}
+u_{ia}g^{ab}R_{bklj}.
\end{split}
\end{equation}
We write $f_z=\fracp{f}{u}$ and $\trw=w^{ij}g_{ij}$.
The letter $c$ denotes estimated positive constants and may change
its value from line to line. It is used so that increasing $c$
keeps the estimates valid. We use $(c_j)$, $(c^k)$, \ldots{} to 
denote estimates tensors.
\end{notation}


\section{Uniform $C^0$-Estimates}
\label{C0 sec}

The techniques of this section are quite standard,
but they simplify the $C^0$-estimates used before for
Schouten tensor equations, see \cite[Prop.\ 3]{ViaclovskyCAG}.
Here, we interpolate between the expressions for the Schouten 
tensors rather than between
the functions inducing the conformal deformations.

For the existence proof
of the regularized problem, we apply a mapping degree argument.
In view of our sub- and supersolutions, we only have to 
ensure that we can apply the maximum principle 
or the Hopf boundary point lemma at a point,
where a solution touches a barrier for the first time during
the deformation associated with the mapping degree argument
to prove $C^0$-estimates. 
Note that $u$ can touch $\ol u$ only in those points from
below where $\ol u$ is admissible. Compare this to
\cite{CGScalar}. Without loss of generality, we may
assume that $u$ touches $\ul u$ from above. Here, touching 
means $u=\ul u$ and $\nabla u=\nabla\ul u$ at a point,
so our considerations include the case of touching at
the boundary.
It suffices to prove an inequality of the form
\begin{equation}\label{ell inequ}
0\le a^{ij}(\ul u-u)_{ij}+b^i(u-\ul u)_i+d(\ul u-u)
\end{equation}
with positive definite $a^{ij}$. 

Define
$$S_{ij}^\psi[v]=v_{ij}+\psi v_iv_j-\tfrac12\psi\vert\nabla
v\vert^2g_{ij}+T_{ij}.$$
We apply the mean value theorem and get for a symmetric 
positive definite tensor $a^{ij}$ and a function $d$
\begin{align*}
0\le&\log\det S^\psi_{ij}[\ul u]-\log\det S^\psi_{ij}[u]
-f(x,\ul u)+f(x,u)\\
=&\int\limits_0^1\frac{d}{dt}\log\det\left\{
tS^\psi_{ij}[\ul u]+(1-t)S^\psi_{ij}[u]\right\}dt
-\int\limits_0^1\frac d{dt}f(x,t\ul u+(1-t)u)dt\\
=&a^{ij}\left(\left(\ul u_{ij}+\psi\ul u_i\ul u_j
-\tfrac12\psi\vert\nabla\ul u\vert^2g_{ij}\right)
-\left(u_{ij}+\psi u_iu_j-\tfrac12\psi
\vert\nabla u\vert^2g_{ij} \right)\right)\\
&+d\cdot(\ul u-u).
\end{align*}
The first integral is well-defined as the set of 
positive definite tensors is convex.
We have $\vert\nabla\ul u\vert^2
-\vert\nabla u\vert^2=\langle\nabla(\ul u-u),\nabla(\ul u+u)\rangle$
and
\begin{align*}
a^{ij}(\ul u_i\ul u_j-u_iu_j)=&a^{ij}\int\limits_0^1\frac d{dt}
\left((t\ul u_i+(1-t)u_i)(t\ul u_j+(1-t)u_j)\right)dt\\
=&2a^{ij}\int\limits_0^1(t\ul u_j+(1-t)u_j)dt
\cdot(\ul u-u)_i,
\end{align*}
so we obtain an inequality of the form \eqref{ell inequ}.
Thus, we may assume in the following that we have 
$\ul u\le u\le\ol u$.


\section{Uniform $C^1$-Estimates}
\label{C1 sec}

\begin{lemma}\label{unif C1}
An admissible solution of \eqref{psi eqn} has uniformly bounded
gradient. 
\end{lemma}
\begin{proof}
We apply a method similar to \cite[Lemma 4.2]{OSMathZ}. Let
$$W=\tfrac12\log\vert\nabla u\vert^2+\mu u$$
for $\mu\gg1$ to be fixed. Assume that W attains its
maximum over $M^n$ at an interior point $x_0$. 
This implies at $x_0$
$$0=W_i=\frac{u^ju_{ji}}{\vert\nabla u\vert^2}+\mu u_i$$
for all $i$. Multiplying with $u^i$ and using 
admissibility gives
\begin{align*}
0=&u^iu^ju_{ij}+\mu\vert\nabla u\vert^4\\
\ge&-\psi\vert\nabla u\vert^4+\tfrac12\psi\vert\nabla u\vert^4
-c\vert\nabla u\vert^2-\lambda\vert\nabla u\vert^2
+\mu\vert\nabla u\vert^4.
\end{align*}
The estimate follows for sufficiently large $\mu$ as
$\lambda$, see \eqref{mod sub sol}, does not depend on $\psi$.
If $W$ attains its maximum at a boundary point $x_0$, we 
introduce normal coordinates such that $W_n$ corresponds
to a derivative in the direction of the inner unit normal.
We obtain in this case $W_i=0$ for $i<n$ and $W_n\le0$
at $x_0$. As the boundary values of $u$ and $\ul u$ 
coincide and $u\ge\ul u$, we may assume that
$u_n\ge0$. Otherwise, $0\ge u_n\ge\ul u_n$ and 
$u_i=\ul u_i$, so a bound for $\vert\nabla u\vert$
follows immediately. Thus we obtain
$0\ge u^iW_i$ and the rest of the proof is identical
to the case where $W$ attains its maximum in the
interior.
\end{proof}
Note that in order to obtain uniform $C^1$-estimates,
we used admissibility, but did not differentiate \eqref{f eqn}.


\section{$C^2$-Estimates}
\label{C2 sec}

\subsection{$C^2$-Estimates at the Boundary}
Boundary estimates for an equation of the form
$\det(u_{ij}+S_{ij})=f(x)$ have been considered
in \cite{CNS1}. It is straight forward to handle 
the additional term that is independent of $u$ in
the determinant and to use subsolutions like in
\cite{BGuanTrans,GuanSpruck,OSMathZ,NehringCrelle}. We want to
point out, that we were only able to obtain 
estimates for the second derivatives of $u$ at
the boundary by introducing $\psi$ and thus removing
gradient terms of $u$ in the determinant near
the boundary. The $C^2$-estimates at the boundary
are very similar to \cite{OSMathZ}. We do not 
repeat the proofs for the double tangential
and double normal estimates, but repeat that for the 
mixed tangential normal derivatives as we can
slightly streamline this part.
Our method does not imply uniform a priori estimates
at the boundary as we look only at small neighborhoods
of the boundary depending on the regularization
or, more precisely, on the set, where $\psi=0$.

\begin{lemma}[Double Tangential Estimates]\label{doub tang est}
An admissible solution of \eqref{psi eqn} has uniformly 
bounded partial second tangential derivatives,
i.\ e.\ for tangential directions $\tau_1$ and $\tau_2$,
$u_{,ij}\tau_1^i\tau_2^j$ is uniformly bounded. 
\end{lemma}
\begin{proof}
This is identical to \cite[Section 5.1]{OSMathZ},
but can also be found at various other places.
It follows directly by differentiating the boundary 
condition twice tangentially.
\end{proof}

\begin{lemma}[Mixed Estimates]
An admissible solution of \eqref{psi eqn} has uniformly
bounded partial second mixed tangential normal derivatives,
i.\ e.\ for a tangential direction $\tau$ and 
for the inner unit normal $\nu$, $u_{,ij}\tau^i\nu^j$
is uniformly bounded.
\end{lemma}
\begin{proof}
The proof is similar to \cite[Section 5.2]{OSMathZ}.
The main differences are as follows.
The modified definition of the linear operator
$T$ in \eqref{TL} clarifies the relation between $T$ 
and the boundary
condition. The term $T_{ij}$ does (in general) 
not vanish in a fixed
boundary point for appropriately chosen coordinates.
In \cite{OSMathZ}, we could choose such coordinates.
Similarly, we choose coordinates such that the
Christoffel symbols become small near a fixed
boundary point. 
Here, we can add and subtract the term $T_{ij}$ in
\eqref{L theta est} as it is independent of $u$.
Finally, we explain here more explicitly
how to apply the inequality for geometric and arithmetic
means in \eqref{big matrix}.

Fix normal coordinates around a point $x_0\in\partial M^n$,
so $g_{ij}(x_0)$ equals the Kronecker delta and the
Christoffel symbols fulfill $\left\vert\Gamma^k_{ij}
\right\vert\le c\dist(\cdot,x_0)=c\vert x-x_0\vert$, where the 
distance is measured in the flat metric 
using our chart, but is equivalent to the distance
with respect to the background metric.
Abbreviate the first $n-1$ coordinates by $\hat x$ and
assume that $M^n$ is locally given by
$\{x^n\ge\omega(\hat x)\}$ for a smooth function $\omega$.
We may assume that $(0,\omega(0))$ corresponds to the
fixed boundary point $x_0$ and $\nabla\omega(0)=0$.
We restrict
our attention to a neighborhood of $x_0$,
$\Omega_\delta=\Omega_\delta(x_0)=M^n\cap
B_\delta(x_0)$ for $\delta>0$ to be fixed sufficiently
small, where $\psi=0$. Thus the equation takes the form 
\begin{equation}\label{bdry eqn}
\log\det(u_{ij}+T_{ij})=\log\det\left(u_{,ij}
-\Gamma^k_{ij}u_k+T_{ij}\right)=f(x,u).
\end{equation}
Assume furthermore that $\delta>0$ is chosen so 
small that the distance function 
to $\partial\Omega$ is smooth in $\Omega_\delta$.

We differentiate the boundary condition tangentially
\begin{equation}\label{bdry d1}
0=(u-\ul u)_{,t}(\hat x, \omega(\hat x))
+(u-\ul u)_{,n}(\hat x, \omega(\hat x))
\omega_{,t}(\hat x),\quad t<n.
\end{equation}
Differentiating \eqref{bdry eqn} yields
\begin{equation}\label{d1a}
w^{ij}\left(u_{,ijk}-\Gamma^l_{ij}u_{,lk}\right)
=f_k+f_zu_k+w^{ij}\left(\Gamma^l_{ij,k}u_l-T_{ij,k}
\right).
\end{equation}
This motivates the definition of the differential
operators $T$ and $L$. Here $t<n$ is fixed and
$\omega$ is evaluated
at the projection of $x$ to the first $n-1$ 
components
\begin{equation}\label{TL}
\begin{split}
Tv:=&v_t+v_n\omega_t,\quad t<n,\\
Lv:=&w^{ij}v_{,ij}-w^{ij}\Gamma^l_{ij}v_l.
\end{split}
\end{equation}
On $\partial M^n$, we have $T(u-\ul u)=0$, so we obtain
$$\vert T(u-\ul u)\vert\le c(\delta)\cdot
\vert x-x_0\vert^2\quad\text{on~}\partial\Omega_\delta.$$
Derivatives of $\ul u$ are a priorily bounded, thus
$$\vert LT(u-\ul u)\vert\le c\cdot\left(1+\trw\right)
\quad{in~}\Omega_\delta.$$
Set $d:=\dist(\cdot,\partial M^n)$, measured in the
Euclidean metric of the fixed coordinates. We 
define for $1\gg\alpha>0$ and $\mu\gg1$ to be chosen
$$\theta:=(u-\ul u)+\alpha d-\mu d^2.$$
The function $\theta$ will be the main part of our
barrier. As $\ul u$ is admissible, there exists
$\epsilon>0$ such that
$$\ul u_{,ij}-\Gamma^l_{ij}\ul u_l+T_{ij}\ge
3\epsilon g_{ij}.$$
We apply the definition of $L$ 
\begin{equation}\label{L theta est}
\begin{split}
L\theta=&w^{ij}\left(u_{,ij}-\Gamma^l_{ij}u_l+T_{ij}\right)
-w^{ij}\left(\ul u_{,ij}-\Gamma^l_{ij}\ul u_l+T_{ij}\right)\\
&+\alpha w^{ij}d_{,ij}-\alpha w^{ij}\Gamma^l_{ij}d_l\\
&-2\mu d w^{ij}d_{,ij}-2\mu w^{ij}d_id_j
+2\mu d w^{ij}\Gamma^l_{ij}d_l
\end{split}
\end{equation}
We have $w^{ij}\left(u_{,ij}-\Gamma^l_{ij}u_l+T_{ij}\right)
=w^{ij}w_{ij}=n$. Due to the admissibility
of $\ul u$, we get $-w^{ij}\left(\ul u_{,ij}-\Gamma^l_{ij}
\ul u_l+T_{ij}\right)\le-3\epsilon\trw$ . 
We fix $\alpha>0$ sufficiently small and obtain
$$\alpha w^{ij}d_{,ij}-\alpha w^{ij}\Gamma^l_{ij}d_l
\le\epsilon\trw.$$
Obviously, we have
$$-2\mu dw^{ij}d_{,ij}+2\mu dw^{ij}\Gamma^l_{ij}d_l
\le c(\mu\delta)\trw.$$
To exploit the term $-2\mu w^{ij}d_id_j$, we use that
$\vert d_i-\delta^n_i\vert\le c\cdot\vert x-x_0\vert\le
c\cdot\delta$, so 
$$-2\mu w^{ij}d_id_j\le-\mu w^{nn}+c(\mu\delta)
\max\limits_{k,\,l}\left\lvert w^{kl}\right\rvert.$$
As $w^{ij}$ is positive definite, we obtain by testing
$\left(\begin{matrix}w^{kk} & w^{kl}\\ w^{kl} & w_{ll}
\end{matrix}\right)$ with the vectors $(1,1)$ and $(1,-1)$
that $\left\lvert w^{kl}\right\rvert\le\trw$. Thus
\eqref{L theta est} implies
\begin{equation}\label{L theta est1}
L\theta\le-2\epsilon\trw-\mu w^{nn}+c+c(\mu\delta)\trw
\end{equation}
We may assume that $\left(w^{ij}\right)_{i,\,j<n}$ 
is diagonal. Then 
\begin{equation}\label{big matrix}
\begin{split}
e^{-f}=\det\left(w^{ij}\right) =&
\det\left(\begin{array}{ccccc}
w^{11} & 0       & \cdots & 0             & w^{1n}   \\
0      & \ddots  & \ddots & \vdots        & \vdots   \\
\vdots & \ddots  & \ddots & 0             & \vdots   \\
0      & \cdots  & 0      & w^{n-1\, n-1} & w^{n-1\, n}\\
w^{1n} & \cdots  & \cdots & w^{n-1\, n}   & w^{nn}   \\
\end{array}\right)\\[.3em]
=& \prod_{i=1}^n w^{ii} \:-\: \sum_{i<n} \left|w^{ni}\right|^2 
\:
\prod_{\genfrac{}{}{0pt}{}{j\neq i}{j<n}} w^{jj} 
\;\leq\; 
\prod_{i=1}^n w^{ii} \;.
\end{split}
\end{equation}
implies that $\trw$ tends to infinity if $w^{nn}$ tends 
to zero. So we can fix $\mu\gg1$ such that the absolute
constant in \eqref{L theta est1}
can be absorbed. Note also that
the geometric arithmetic means inequality implies 
$$\tfrac1n\trw=\tfrac1n\sum\limits_{i=1}^n 
w^{ii}\ge\left(\prod\limits_{i=1}^n w^{ii}\right)^{1/n},$$
so \eqref{big matrix} yields a positive lower bound
for $\trw$. Finally, we
fix $\delta=\delta(\mu)$ sufficiently small and use 
\eqref{L theta est1} to deduce that
\begin{equation}\label{L theta fin}
L\theta\le-\epsilon\trw.
\end{equation}
We may assume that $\delta$
is fixed so small that $\theta\ge0$ in $\Omega_{\delta}$.

Define for $A,\,B\gg1$ the function
$$\Theta^{\pm}:=A\theta+B\vert x-x_0\vert^2
\pm T(u-\ul u).$$
Our estimates imply that $\Theta^{\pm}\ge0$ on 
$\partial\Omega_\delta$ for $B\gg1$ fixed sufficiently
large and $L\Theta^{\pm}\le0$ in $\Omega_\delta$, when
$A\gg1$, depending also on $B$, is fixed sufficiently
large. Thus the maximum principle implies that
$\Theta^{\pm}\ge0$ in $\Omega_\delta$. As $\Theta^{\pm}
(x_0)=0$, we deduce that $\Theta^{\pm}_{,n}\ge0$, so we
obtain a bound for $(Tu)_{,n}$ and the lemma follows.
\end{proof}

\begin{lemma}[Double Normal Estimates]
An admissible solution of \eqref{psi eqn} has uniformly
bounded partial second normal derivatives,
i.\ e.\ for the inner unit normal $\nu$, 
$u_{,ij}\nu^i\nu^j$ is uniformly bounded.
\end{lemma}
\begin{proof}
The proof is identical to \cite[Section 5.3]{OSMathZ}.
Note however, that the notation there is slightly
different. There $-u_{,ij}+a_{ij}$ is positive
definite instead of $u_{,ij}-\Gamma^k_{ij}u_k+T_{ij}$
here. 
\end{proof}

\subsection{Interior $C^2$-Estimates}

\begin{lemma}[Interior Estimates]\label{int C2}
An admissible solution of \eqref{psi eqn} has uniformly
bounded second derivatives.
\end{lemma}
\begin{proof}
Note the admissibility implies that $w_{ij}$ is positive
definite. This implies a lower bound on the eigenvalues
of $u_{ij}$. 

For $\lambda\gg1$ to be chosen sufficiently large, 
we maximize the functional 
$$W=\log\left(w_{ij}\eta^i\eta^j\right)+\lambda\chi$$
over $M^n$ and all $\left(\eta^i\right)$ with 
$g_{ij}\eta^i\eta^j=1$. In view of the boundary estimates
obtained above, we may assume that $W$ attains its maximum
at an interior point $x_0$ of $M^n$.
As in \cite{CGJDG1996} we may choose normal coordinates 
around $x_0$ and an appropriate extension of 
$\left(\eta^i\right)$ corresponding to the maximum
value of $W$. In this way, we can pretend that
$w_{11}$ is a scalar function that equals $w_{ij}\eta^i
\eta^j$ at $x_0$ and we obtain
\begin{align}\label{Wi}
0=W_i=&\frac1{w_{11}}w_{11;i}+\lambda\chi_i\\
\intertext{and}
0\ge W_{ij}=&\frac1{w_{11}}w_{11;ij}-\frac1{w_{11}^2}
w_{11;i}w_{11;j}+\lambda\chi_{ij}\label{Wij}
\end{align}
in the matrix sense, $1\le i,\,j\le n$. Here and below,
all quantities are evaluated at $x_0$. We may 
assume that $w_{ij}$ is diagonal and $w_{11}\ge 1$.
Differentiating \eqref{psi eqn} yields
\begin{align}\label{d1}
w^{ij}w_{ij;k}=&f_k+f_z u_k,\\
w^{ij}w_{ij;11}-w^{ik}w^{jl}w_{ij;1}w_{kl;1}=&
f_{11}+2f_{1z}u_1+f_{zz}u_1u_1+f_zu_{11}.\label{d2}
\end{align}
Combining the convexity assumption on $\chi$,
\eqref{Wij} and \eqref{d2} gives
\begin{equation}\label{wijWij}
\begin{split}
0\ge&\frac1{w_{11}}w^{ij}w_{11;ij}
-\frac1{w_{11}^2}w^{ij}w_{11;i}w_{11;j}+\lambda\trw\\
=&\frac1{w_{11}}w^{ij}(w_{11;ij}-w_{ij;11})\\
&+\frac1{w_{11}}w^{ik}w^{jl}w_{ij;1}w_{kl;1}
-\frac1{w_{11}^2}w^{ij}w_{11;i}w_{11;j}\\
&+\frac1{w_{11}}(f_{11}+2f_{1z}u_1+f_{zz}u_1u_1+f_zu_{11})
+\lambda\trw.
\end{split}
\end{equation}
It will be convenient to decompose $w_{ij}$ as follows
\begin{equation}\label{r intro}
\begin{split}
w_{ij}=&u_{ij}+r_{ij},\\
r_{ij}=&\psi u_iu_j-\tfrac12\psi\abs{\nabla u}^2g_{ij}+T_{ij}.
\end{split}
\end{equation}
The quantity $r_{ij}$ is a priorily bounded, so 
the right-hand side of \eqref{d2} is bounded from
below by $-c(1+w_{11})$. 

Let us first consider some terms involving at most third
derivatives of $u$
\begin{equation}\label{le three}
\begin{split}
w^{ik}w^{jl}&w_{ij;1}w_{kl;1}-\frac1{w_{11}}
w^{ij}w_{11;i}w_{11;j}
\ge\frac1{w_{11}}w^{ij}(w_{i1;1}w_{j1;1}-w_{11;i}w_{11;j})\\
=&\frac1{w_{11}}w^{ij}((u_{i11}+r_{i1;1})(u_{j11}+r_{j1;1})
-(u_{11i}+r_{11;i})(u_{11j}+r_{11;j}))\\
\ge&\frac1{w_{11}}w^{ij}(u_{i11}u_{j11}-u_{11i}u_{11j}
+2u_{i11}r_{j1;1}-2u_{11i}r_{11;j}-r_{11;i}r_{11;j}).\\
\end{split}
\end{equation}
We will bound each term on the right-hand side individually.
The term $r_{11;i}$ is of the form $c_i+c^ku_{ki}$. 
We rewrite $u_{ki}=w_{ki}-r_{ki}$,
use $w^{ij}w_{jk}=\delta^i_k$ and obtain 
$$\left\lvert\frac1{w_{11}}w^{ij}r_{11;i}r_{11;j}\right\rvert
\le\frac1{w_{11}}c\left(1+w_{11}+\trw\right).$$
Note that $w_{ij}$ is diagonal, so the maximality of
$W$ implies $\abs{w_{ij}}\le w_{11}$ for any $i,\,j$.
We use \eqref{r intro}, \eqref{Wi} and rewrite
$r_{11;i}$ as above
\begin{align*}
-2\frac1{w_{11}}w^{ij}u_{11i}r_{11;j}=&-2\frac1{w_{11}}w^{ij}
(w_{11i}-r_{11;i})r_{11;j}\\
=&2\lambda w^{ij}\chi_ir_{11;j}+2\frac1{w_{11}}w^{ij}
r_{11;i}r_{11;j}\\
\ge&-c\lambda\left(1+\trw\right)
-\frac1{w_{11}}c\left(1+w_{11}+\trw\right).
\end{align*}
To estimate the next term, we use \eqref{interchange},
\eqref{r intro}, \eqref{Wi} and the fact that the second
derivatives of $u$ in $r_{j1;1}$ appear with a factor $\psi$
\begin{align*}
\frac2{w_{11}}w^{ij}u_{i11}r_{j1;1}=&
\frac2{w_{11}}w^{ij}\left(w_{11;i}-r_{11;i}+u_ag^{ab}
R_{b1i1}\right)r_{j1;1}\\
=&-2\lambda w^{ij}\chi_i\left(c_j+\psi c^k_jw_{k1}+c^kw_{kj}\right)\\
&-\frac c{w_{11}}w^{ij}\left(c_i+c^kw_{ki}\right)\left(c_j+
\psi c^k_jw_{k1}+c^kw_{kj}\right)\\
\ge&-c\lambda\left(1+\trw+\psi w_{11}\trw\right)
-c\left(1+\trw\right).
\end{align*}
We interchange third covariant derivatives and get
\begin{align*}
\frac1{w_{11}}w^{ij}&(u_{i11}u_{j11}-u_{11i}u_{11j})\displaybreak[1]\\
=&\frac1{w_{11}}w^{ij}\left(u_{i11}u_{j11}-
\left(u_{i11}+u_ag^{ab}R_{b11i}\right)\left(u_{j11}+u_cg^{cd}
R_{d11j}\right)\right)\displaybreak[1]\\
\ge&-2\frac1{w_{11}}w^{ij}u_{i11}u_ag^{ab}R_{b11j}
-c\frac1{w_{11}}\trw\displaybreak[1]\\
=&2\lambda w^{ij}\chi_iu_ag^{ab}R_{b11j}
+2\frac1{w_{11}}w^{ij}r_{i1;1}u_ag^{ab}R_{b11j}
-c\frac1{w_{11}}\trw\displaybreak[1]\\
\ge&-c(1+\lambda)\left(1+\trw\right).
\end{align*}
Recall that $\trw$ is bounded below by a positive 
constant.
We employ \eqref{le three} and get the estimate
\begin{equation}\label{three est}
\frac1{w_{11}}w^{ik}w^{jl}w_{ij;1}w_{kl;1}
-\frac1{w_{11}^2}w^{ij}w_{11;i}w_{11;j}\ge
-c\left(1+\lambda\psi+\frac{\lambda}{w_{11}}\right)\trw.
\end{equation}

Next, we consider the terms in \eqref{wijWij} involving
fourth derivatives. Equation \eqref{interchange} implies
\begin{align*}
u_{11ij}=&u_{ij11}+u_{a1}g^{ab}R_{bi1j}+u_ag^{ab}R_{bi1j;1}
+u_{1a}g^{ab}R_{bij1}+u_{ia}g^{ab}R_{b1j1}\\
&+u_{aj}g^{ab}R_{b11i}+u_ag^{ab}R_{b11i;j}\\
\ge&u_{ij11}-c_{ij}(1+w_{11}).
\end{align*}
We use \eqref{r intro}
\begin{align*}
w^{ij}(w_{11;ij}&-w_{ij;11})=w^{ij}(u_{11ij}-u_{ij11})
+w^{ij}(r_{11;ij}-r_{ij;11})\displaybreak[1]\\
\ge&w^{ij}(r_{11;ij}-r_{ij;11})-c w_{11}\trw\\
=&w^{ij}\left(\psi_{ij}u_1^2+4\psi_iu_1u_{1j}+2\psi u_{1j}u_{1i}
+2\psi u_1u_{1ij}\right)\\
&+w^{ij}\left(-\psi_{11}u_iu_j-4\psi_1u_{i1}u_j-2\psi u_{1i}u_{1j}
-2\psi u_iu_{j11}\right)\\
&+w^{ij}\left(-\tfrac12\psi_{ij}\abs{\nabla u}^2g_{11}
-2\psi_iu^ku_{kj}g_{11}-\psi u^k_ju_{ki}g_{11}
-\psi u^ku_{kij}g_{11}\right)\\
&+w^{ij}\left(\tfrac12\psi_{11}\abs{\nabla u}^2g_{ij}
+2\psi_1u^ku_{k1}g_{ij}+\psi u^k_1u_{k1}g_{ij}
+\psi u^ku_{k11}g_{ij}\right)\\
&+w^{ij}\left(T_{11;ij}-T_{ij;11}\right)-c w_{11}\trw.
\end{align*}
Some terms cancel. We use \eqref{r intro} and the fact
that $w^{ij}$ is the inverse of $w_{ij}$. Then we interchange
covariant third derivatives \eqref{interchange} and
employ once again \eqref{r intro}
\begin{align*}
w^{ij}(w_{11;ij}&-w_{ij;11})\ge\\
\ge&w^{ij}\left(2\psi u_1u_{1ij}-2\psi u_iu_{j11}
-\psi u^ku_{kij}g_{11}+\psi u^ku_{k11}g_{ij}\right)\\
&+w^{ij}\left(-\psi u^k_ju_{ki}g_{11}
+\psi u^k_1u_{k1}g_{ij}\right)-cw_{11}\trw\displaybreak[1]\\
=&2\psi u_1w^{ij}u_{ij1}+2\psi u_1w^{ij}u_ag^{ab}R_{bi1j}\\
&-\psi g_{11}u^kw^{ij}u_{ijk}-\psi g_{11}u^kw^{ij}u_ag^{ab}R_{bikj}\\
&-2\psi u_iw^{ij}u_{11j}-2\psi u_iw^{ij}u_ag^{ab}R_{b1j1}\\
&+\psi u^ku_{11k}\trw+\psi u^ku_ag^{ab}R_{b1k1}\trw\\
&-\psi g_{11}w^{ij}(w_{ik}-r_{ik})(w_{jl}-r_{jl})g^{kl}\\
&+\psi(w_{1k}-r_{1k})(w_{1l}-r_{1l})g^{kl}\trw
-c w_{11}\trw.\\
\end{align*}
We replace third derivatives of $u$ by derivatives of $w_{ij}$.
Equations \eqref{d1} and \eqref{Wi} allow to replace these
terms by terms involving at most second derivatives of $u$
\begin{align*}
w^{ij}(w_{11;ij}&-w_{ij;11})\ge\displaybreak[1]\\
\ge&2\psi u_1w^{ij}w_{ij;1}-2\psi u_1w^{ij}r_{ij;1}
-\psi g_{11}u^kw^{ij}w_{ij;k}+\psi g_{11}u^kw^{ij}r_{ij;k}\\
&-2\psi u_iw^{ij}w_{11;j}+2\psi u_iw^{ij}r_{11;j}
+\psi u^kw_{11;k}\trw-\psi u^kr_{11;k}\trw\\
&+\psi w_{11}^2\trw-c w_{11}\trw\displaybreak[1]\\
\ge&2\lambda\psi w_{11}w^{ij}u_i\chi_j-\lambda\psi w_{11}
u^k\chi_k\trw+\psi w_{11}^2\trw-c w_{11}\trw\displaybreak[1]\\
\ge&-c\lambda\psi w_{11}\trw+\psi w_{11}^2\trw
-cw_{11}\trw.
\end{align*}
This gives
\begin{equation}\label{four est}
\frac1{w_{11}}w^{ij}(w_{11;ij}-w_{ij;11})\ge
-c\lambda\psi\trw+\psi w_{11}\trw-c\trw.
\end{equation}
We estimate the respective terms in \eqref{wijWij}
using \eqref{three est} and \eqref{four est} and obtain
\begin{equation}\label{together}
0\ge\left\{\psi(w_{11}-c\lambda)
+\left(\lambda-c-\frac{c\lambda}{w_{11}}\right)
\right\}\trw.
\end{equation}
Assume that all $c$'s in \eqref{together} are equal.
Now we fix $\lambda$ equal to $c+1$. 
Then \eqref{together} implies
that $w_{11}$ is bounded above.
\end{proof}

\bibliographystyle{/usr/people/schnuere/os/tools/amsplain}
\bibliography{/usr/people/schnuere/os/tools/biblio}

\end{document}